\documentclass[leqno]{amsart}
\usepackage{amsfonts,amssymb,amsmath,amsgen,amsthm}
\usepackage{hyperref}
\usepackage{color}
\newcommand{\msc}[2][2000]{%
  \let\@oldtitle\@title%
  \gdef\@title{\@oldtitle\footnotetext{#1 \emph{Mathematics subject
        classification.} #2}}%
}

\theoremstyle{plain}
\newtheorem{theorem}{Theorem} [section]

\newtheorem{assumption}[theorem]{Assumption}
\newtheorem{lemma}[theorem]{Lemma}

\newtheorem{proposition}[theorem]{Proposition}

\theoremstyle{remark}
\newtheorem{remark}[theorem]{Remark}

\def\C{{\mathbb C}}
\def\R{{\mathbb R}}
\def\N{{\mathbb N}}
\def\Sch{{\mathcal S}}
\def\O{\mathcal O}

\def\({\left(}
\def\){\right)}
\def\<{\left\langle}
\def\>{\right\rangle}
\def\le{\leqslant}
\def\ge{\geqslant}

\def\d{{\partial}}
\def\eps{\varepsilon}

\def\si{{\sigma}}

\def\eik{\phi_{\rm eik}}

\DeclareMathOperator{\RE}{Re}
\DeclareMathOperator{\IM}{Im}
\DeclareMathOperator{\diver}{div}
\DeclareMathOperator{\supp}{supp}

\numberwithin{equation}{section}

\begin{document}
\title[Multiphase WKB analysis for NLS]{On nonlinear effects in
  multiphase WKB analysis for the 
  nonlinear Schr\"odinger equation}  
\author[R. Carles]{R\'emi Carles}
\address{Univ Rennes, CNRS\\ IRMAR - UMR 6625\\ F-35000
  Rennes, France}
\email{Remi.Carles@math.cnrs.fr}

\begin{abstract}
  We consider the Schr\"odinger equation with an external potential
  and a cubic nonlinearity, in the semiclassical limit. The initial
  data are sums of WKB states, with smooth phases and smooth, compactly
  supported initial amplitudes, with disjoint supports. We show that
  like in the linear case, a superposition principle holds on some time
  interval independent of the semiclassical parameter, in several
  r\'egimes in term of 
  the size of initial data with respect to the semiclassical parameter. 
 For large data, we invoke properties of compressible Euler
  equations. For smaller data, we show that there may be no
 nonlinear  interferences on some time interval independent of the
 semiclassical parameter, and interferences for later time, thanks to
 explicit computations available for particular phases. 
\end{abstract}

\thanks{This work was supported by Centre Henri Lebesgue,
program ANR-11-LABX-0020-0. A CC-BY public copyright license has
  been applied by the author to the present document and will be
  applied to all subsequent versions up to the Author Accepted
  Manuscript arising from this submission.}

\maketitle

\section{Introduction}
\label{sec:intro}
\subsection{Setting}
\label{sec:setting}

We consider the cubic defocusing Schr\"odinger equation on $\R^d$, $d\ge
1$, in the semiclassical r\'egime
\begin{equation*}
  i\eps\d_t \psi^\eps +\frac{\eps^2}{2}\Delta
  \psi^\eps=V\psi^\eps+|\psi^\eps|^2 \psi^\eps. 
\end{equation*}
The potential $V=V(x)$ is supposed real-valued, smooth, and at most
quadratic:
\begin{equation}
  \label{eq:V}
  V\in C^\infty(\R^d,\R),\quad \d^\alpha V\in L^\infty(\R^d),\quad
  \forall \alpha\in \N^d, \ |\alpha|\ge 2. 
\end{equation}
Typical examples are $V=0$, $V$ linear ($V(x)=E\cdot x$), $V$ harmonic
($V(x)=\omega^2|x|^2$), $V\in \Sch(\R^d)$, or any sum of such potentials. 
As initial data, we consider the sum of  WKB states:
 \begin{equation*}
 \psi^\eps(0,x)=\eps^{\gamma/2}\sum_{j=1}^N
 \alpha_j(x)e^{i\varphi_j(x)/\eps} , 
\end{equation*}
with $\gamma \ge 0$. The value of $\gamma$ measures the size of the
initial data, and thus the importance of nonlinear effects in the
semiclassical limit $\eps\to 0$. The case $\gamma =0$ is supercritical in terms of WKB
analysis:  the evolution
of the phase describing the rapid oscillation is given by an eikonal
equation which involves the leading order amplitude, and a standard
application of the WKB asymptotic expansion leads to systems which are
not closed, no matter how many correcting terms are considered  (see
e.g. \cite[Chapter~1]{CaBook2} or \cite{PGX93}). 
\smallbreak

In order to consider solutions of size $\O(1)$ in terms of $\eps$, we
consider the unknown function $u^\eps = \eps^{-\gamma/2} \psi^\eps$
instead of $\psi^\eps$, where $u^\eps$ thus solves
\begin{equation}
  \label{eq:NLS}
  i\eps\d_t u^\eps +\frac{\eps^2}{2}\Delta
  u^\eps=Vu^\eps+\eps^\gamma|u^\eps|^2 u^\eps,
\end{equation}
with initial data 
 \begin{equation}\label{eq:CI}
 u^\eps(0,x)=  u_0^\eps(x):=\sum_{j=1}^N
 \alpha_j(x)e^{i\varphi_j(x)/\eps} .
\end{equation}
The assumptions regarding the initial data are listed below:
\begin{assumption}\label{hyp:CI}
  The phases are smooth and real-valued, $\varphi_j\in
C^\infty(\R^d,\R)$. The initial amplitudes
are smooth and compactly supported: $\alpha_j\in
C_0^\infty(\R^d,\C)$, with pairwise disjoint supports,
\begin{equation*}
  \supp \alpha_{j_1} \cap \supp \alpha_{j_2}=\emptyset,\quad j_1\not
  =j_2. 
\end{equation*}
\end{assumption}
The case $N=1$, referred to as \emph{monokinetic case}, is well
understood for \emph{short time}, as we recall below, in the sense
that the asymptotic behavior of $u^\eps$ as $\eps\to
0$ is  described precisely, locally in time on some interval
independent of 
$\eps$. The large time behavior is, in general, unknown; the
one-dimensional case, with $V=0$, is an exception, since it is completely
integrable, see e.g. \cite{JLM,Miller2016}. Consider the case
$\gamma=0$. When
$V\equiv 0$, the leading
order asymptotic description involves the compressible Euler equation
\begin{equation}
  \label{eq:euler}
  \left\{
    \begin{aligned}
     & \d_t \rho+\diver \(\rho v\)=0,\\
    & \d_t v + v\cdot \nabla v+\nabla \rho=0.
          \end{aligned}
 \right.
\end{equation}
This equation is quasilinear, while \eqref{eq:NLS} is semilinear (the
nonlinear term is viewed as a perturbation when solving the Cauchy
problem). In Section~\ref{sec:mono}, we recall how to justify, in this
case, the existence of a WKB approximation of the form
\begin{equation*}
  u^\eps(t,x) =\( {\tt a}(t,x)+\eps {\tt a}_1(t,x)+\ldots +\eps^k{\tt
    a}_k(t,x)\)e^{i\phi(t,x)/\eps}+\O\(\eps^{k+1}\),
\end{equation*}
 in $L^\infty([0,T],L^2\cap L^\infty(\R^d))$, for all $k\ge 0$,
 for some $T>0$ independent of $\eps$. We choose to measure errors in
 $L^2\cap L^\infty$ in the spatial norm, in order to avoid to
 introduce $\eps$-dependent norms  when
 derivatives are involved, due to rapid oscillations. This time $T$
 can be taken as 
 the lifespan of the smooth solution to the Euler equation
 \eqref{eq:euler} with suitable initial data. 
 When $N\ge 2$, the new
question arising is the nonlinear interaction of the WKB
states. As the problem is supercritical, even a formal computation is
a delicate issue: if we plug an approximate solution of the form
\begin{equation*}
  u^\eps_{\rm app}(t,x) = \sum_{j=1}^M  b_j(t,x)e^{i\phi_j(t,x)/\eps} 
\end{equation*}
into \eqref{eq:NLS}, how do we choose $M$ (possibly infinite), and
which equations must be satisfied by the amplitudes $b_j$ and the
phases $\phi_j$?  Surprisingly enough, it turns out that as long as the
solutions of the Euler equations, involved in the description of each
individual initial WKB state, are smooth, there is no interaction, at
arbitrary order in terms of powers of $\eps$. 

\begin{remark}[Infinitely many states] The case
$N=\infty$ may also be addressed, under suitable assumptions on the
growth in space of the phases $\phi_j$ compared to the size of the
support of $\alpha_j$, as $j\to \infty$. More precisely, as will be
clear from the proof of the main result, we can consider the case
$N=\infty$ provided that we may find cutoff
functions $\chi_j$ so that 
\begin{equation*}
  \phi_0=\sum_{j=1}^\infty\varphi_j\chi_j\in H^\infty(\R^d):=\cap_{s>0} H^s(\R^d),
\end{equation*}
or at least in a weaker form if $\phi_0\in H^s(\R^d)$ for some $s>2+d/2$. Another
constraint, in this case, is that we have to find a common lower bound for
the lifespan of all the  approximate solutions $(\phi_j,a_j)$
considered below, an aspect which is obvious
when $N$ is finite, since we consider the minimum of a finite set.  
\end{remark}

\subsection{Main results}

The nonlinear evolution of each initial WKB state will play a crucial role: 
\begin{equation}
 \label{eq:NLSj}
  i\eps\d_t u^\eps_j +\frac{\eps^2}{2}\Delta
  u_j^\eps=Vu_j^\eps+\eps^\gamma|u_j^\eps|^2 u_j^\eps\quad ;\quad u_{j\mid t=0}
  =\alpha_j e^{i\varphi_j/\eps}. 
\end{equation}
Under our assumptions, for fixed initial data, we know that:
\begin{itemize}
\item If $d\le 3$, the equation is energy-subcritical, and for fixed
  $\eps>0$, there exists a unique solution $u^\eps\in
  L^\infty(\R;H^1(\R^d))$, and it is smooth. See e.g. \cite{CazCourant}. 
 \item If $d=4$, the equation is energy-critical: the above conclusion
   is known to 
   remain  when $V=0$ (\cite{RV07}), when $V$ is an isotropic
   quadratic potential (\cite{KVZ09}), or when $V$ is harmonic at
   infinity (\cite{Casey2018}). 
 \item If $d\ge 5$, the equation is energy-supercritical: only a local in time smooth solution is known to
   exist by classical theory.
\end{itemize}
In the cases where the global existence of a smooth solution is not
known, the local existence time might go to zero as $\eps\to 0$, so
the existence of a smooth solution on a time interval independent of
$\eps>0$ is already a nontrivial step. 
The description of the solutions $u_j^\eps$ as $\eps\to 0$ on some
time interval $[0,T_j]$ independent of $\eps$ was evoked above, and is
recalled 
in Sections~\ref{sec:mono} (case $V=0$) and \ref{sec:monopot} ($V$
satisfying \eqref{eq:V}).
Our main result is the following nonlinear superposition principle:

\begin{theorem}\label{theo:superp}
  Let $d\ge 1$, $V$ satisfying \eqref{eq:V}, $\gamma\ge0$, and initial data
  satisfying Assumption~\ref{hyp:CI}. There exists $T^*>0$ independent
  of $\eps\in ]0,1]$ such that \eqref{eq:NLS}-\eqref{eq:CI} has a
  unique solution $u^\eps\in C([0,T^*],H^\infty(\R^d))$. In addition,
  \begin{equation*}
 \sup_{t\in [0,T^*]} \left\|u^\eps(t)-\sum_{j=1}^N
   u_j^\eps(t)\right\|_{L^2\cap L^\infty}=\O\(\eps^k\),\quad \forall k>0,
\end{equation*}
where $u_j^\eps$ is the solution of \eqref{eq:NLSj}. 
\end{theorem}
Let us discuss this result in the supercritical case $\gamma=0$, as it
is the case where Theorem~\ref{theo:superp} may be the most surprising. 
The result follows from a detailed WKB analysis, as
well as a property of finite speed of propagation for the
compressible Euler equation, discovered initially in \cite{MUK86}. The
key feature of our setting  is the compact, disjoint supports of
the initial amplitudes $\alpha_j$. In the case $V=0$, as long as WKB
analysis is valid for each $u_j^\eps$ in \eqref{eq:NLSj}, $u_j^\eps$
remains supported in (essentially) $\supp \alpha_j$ up to
$\O(\eps^\infty)$: all the amplitude terms of the WKB expansion  (at
leading order, as well as correctors at arbitrary order) remain
compactly supported, and amplitudes associated with $u_{j_1}^\eps$ and
$u_{j_2}^\eps$, respectively, with $j_1\not = j_2$, do not interact.  In the
case $V\not\equiv 0$, $u_j^\eps$
remains supported in $\supp \alpha_j$ evolving according to the
classical flow generated by $V$, up to $\O(\eps^\infty)$. In other
words, we recover the same phenomenon, regarding the evolution of supports, as in the linear case (see
e.g. \cite{MaslovFedoryuk,Robert}), even 
though the r\'egime associated to \eqref{eq:NLS} is strongly
nonlinear (of course Theorem~\ref{theo:superp} is trivial in the
linear case, as $u^\eps\equiv \sum u_j^\eps$).  In particular, the initial modes cannot interact at a
``visible'' order before WKB analysis for at least one of the
$u_j^\eps$'s ceases to be valid, that is, before the solution of the
corresponding Euler equation \eqref{eq:euler} breaks down (see however
Section~\ref{sec:chi} for a discussion on the influence of our
proof strategy on this statement). Recent
progress on this 
precise question, \cite{MR4445442,MR4445443,BCLGS-p} (see also \cite{MRRS22} for a relation with the
nonlinear Schr\"odinger equation), suggests that the expected scenario
is rather that of an implosion: the conclusion of
Theorem~\ref{theo:superp} might remain valid even after WKB has ceased
to be valid. 

\begin{remark}[Wigner measures]\label{rem:wigner1}
 Since the proof of Theorem~\ref{theo:superp} relies on WKB analysis,
 it also implies the characterization of Wigner measures. Recall that
 the Wigner transform of $u^\eps$ is defined by 
\begin{equation*} 
  w^\eps(t,x,\xi )
= (2\pi )^{-d}\int_{\R^d} u^\eps\(t,x-\eps\frac{\eta}{2}\)
\overline  u^\eps\(t,x+\eps \frac{\eta}{2}\)\,
e^{i \eta \cdot \xi} d\eta .
\end{equation*}
The position and current densities can be recovered from $w^\eps$, by 
\begin{align*}
  &|u^\eps(t,x)|^2 = \int_{\R^d} w^\eps(t,x,\xi)d\xi,\\
 & \IM\(\eps \overline u^\eps\nabla u^\eps\)(t,x) =
 \int_{\R^d} \xi w^\eps(t,x,\xi)d\xi.  
\end{align*}
A measure $\mu$ is a Wigner
measure associated to $u^\eps$ (there is no uniqueness in general)
if, up to extracting a subsequence, $w^\eps$ converges to $\mu$ as $\eps\to
0$ (see e.g. \cite{GMMP,LionsPaul}). In the context of
Theorem~\ref{theo:superp}, each wave function 
$u_j^\eps$ has a unique Wigner measure, and the sum of these Wigner
measures is the Wigner measure of $u^\eps$. For instance, if $V=\gamma=0$, 
\begin{equation*}
  \mu (t,dx,d\xi) = \sum_{j=1}^N \rho_j(t,x)dx\otimes \delta_{\xi = v_j(t,x)},
\end{equation*}
where $(\rho_j,v_j)$ solves \eqref{eq:euler} with initial data
$(\rho_j,v_j)_{\mid t=0}= (|\alpha_j|^2,\nabla (\chi_j\varphi_j))$,
and $\chi_j\in C_0^\infty(\R^d,[0,1]))$ is (any function) such that
$\chi_j\equiv 1$ on the support of $\alpha_j$. See
Section~\ref{sec:chi} for the dependence of this statement
upon $\chi_j$. 
\end{remark}

\smallbreak

The next result shows that in the weakly nonlinear case $\gamma=1$,
some explicit information is available, in the sense that indeed,
nonlinear interferences are negligible on some time interval $[0,T_0]$
with $T_0>0$ independent of $\eps$, while nonlinear interferences
occur later. 

\begin{proposition}\label{prop:wnlgo}
  Let $d\ge 1$. There exist $k_1,k_2\in \R^d$, and
  $\alpha_1,\alpha_2\in C_0^\infty(\R^d)$ with disjoint
  supports, such that the following holds.  There exist $T_1>0$ 
  and $T_0\in (0,T_1)$ independent of $\eps$, such that the
  solution to
   \begin{equation}\label{eq:wnlgo}
     i\eps\d_t u^\eps +\frac{\eps^2}{2}\Delta u^\eps=\eps |u^\eps|^2
     u^\eps\quad ;\quad u_0^\eps(x)=
     \alpha_1(x)e^{ik_1\cdot x/\eps} + \alpha_2(x)e^{ik_2\cdot x/\eps} ,
   \end{equation}
   satisfies
  \begin{equation*}
    \sup_{t\in [0,T_0]} \left\|u^\eps(t)-\sum_{j=1}^2
   u_j^\eps(t)\right\|_{L^p}=\O\(\eps^k\),\quad \forall k>0,\quad \forall p\in [2,\infty],
\end{equation*}
and
   \begin{equation*}
     \liminf_{\eps\to 0}\sup_{t\in [0,T_1]} \left\|u^\eps(t)-\sum_{j=1}^2
   u_j^\eps(t)\right\|_{L^p}>0,\quad \forall p\in [2,\infty].
\end{equation*}
 \end{proposition}
 The proof of
 Proposition~\ref{prop:wnlgo} relies on explicit computations
 available in this weakly nonlinear case, and the fact that for linear
 oscillations, no caustic appears in the case of a single WKB state:
 the nature of nonlinear interferences is shown in
 Section~\ref{sec:wnlgo}, and consists of nonlinear phase
 modulations. 
In an appendix, we give an alternative argument illustrating another
type of nonlinear
interferences at leading order, consisting of the creation of a new
mode (when $d\ge 2$): starting from three WKB states, a fourth one,
associated with a new phase, may appear by resonant interaction. 
 
 \subsection{Content}
\label{sec:content}

In Section~\ref{sec:mono}, we recall the WKB construction introduced in
\cite{Grenier98} for the case $\gamma=0$, and emphasize the finite
speed of propagation which 
appears in our framework. In Section~\ref{sec:monopot}, we explain how
to adapt the previous approach to the case where $V$ satisfies
\eqref{eq:V}, and address the case $\gamma>0$. In
Section~\ref{sec:separation}, we complete the proof of
Theorem~\ref{theo:superp}. Section~\ref{sec:chi} clarifies the role of
the cutoff functions used in the proof of
Theorem~\ref{theo:superp}. Propositions~\ref{prop:wnlgo} is established
in Section~\ref{sec:wnlgo}. In an appendix, we propose an alternative
proof of Proposition~\ref{prop:wnlgo}, in the case
$d\ge 2$ with $N=3$, showing that there are several sorts of nonlinear
interferences in the weakly nonlinear case.

\section{The monokinetic case without potential}
\label{sec:mono}
In this section, we consider \eqref{eq:NLS}-\eqref{eq:CI} in the
monokinetic $N=1$, and in the supercritical case $\gamma=0$, with
slightly different notations for future reference:
\begin{equation}
  \label{eq:mono}
  i\eps\d_t u^\eps+\frac{\eps^2}{2}\Delta u^\eps=
  |u^\eps|^2u^\eps\quad ;\quad u^\eps_{\mid t=0}=a_0e^{i\phi_0/\eps}.
\end{equation}
In view of the setting of this paper, we assume $a_0,\phi_0\in
C_0^\infty(\R^d)$. In particular, $a_0,\phi_0\in H^\infty(\R^d)$. 
We first consider the case $V\equiv 0$, then introduce the main ideas
that make it possible to incorporate a subquadratic potential $V$.

We recall the main steps to the construction introduced in
\cite{Grenier98} (see also \cite[Section~4.2]{CaBook2}). 
The idea introduced in \cite{Grenier98} consists in writing the
solution to \eqref{eq:mono} as
\begin{equation}\label{eq:reprGrenier}
  u^\eps(t,x) = a^\eps(t,x)e^{i\phi^\eps(t,x)/\eps},
\end{equation}
with $a^\eps$ complex-valued and $\phi^\eps$ real-valued, solving
\begin{equation}\label{eq:Grenier}
  \left\{
    \begin{aligned}
     &\d_t \phi^\eps+\frac{1}{2}|\nabla
     \phi^\eps|^2+|a^\eps|^2=0,\quad \phi^\eps_{\mid t=0} =\phi_0,\\
&\d_t a^\eps+\nabla\phi^\eps\cdot \nabla
a^\eps+\frac{1}{2}a^\eps\Delta \phi^\eps = i\frac{\eps}{2}\Delta
a^\eps,\quad a^\eps_{\mid t=0}=a_0.
    \end{aligned}
\right.
\end{equation}
The key remark is that this leads to a symmetric hyperbolic system,
perturbed be a skew-symmetric term. The hyperbolic system appears when
considering the unknown
\begin{equation*}
  U^\eps =
  \begin{pmatrix}
    \RE a^\eps\\
    \IM a^\eps\\
    \nabla \phi^\eps
  \end{pmatrix}= \begin{pmatrix}
    \RE a^\eps\\
    \IM a^\eps\\
    v^\eps
  \end{pmatrix}.
\end{equation*}
Considering the gradient of the first equation in \eqref{eq:Grenier},
the system can be written
\begin{equation}\label{eq:symL}
  \d_t U^\eps+\sum_{j=1}^d A_j(U^\eps)\d_j U^\eps = \eps LU^\eps,
\end{equation}
with
\begin{equation*}
A(U^\eps,\xi)=\sum_{j=1}^d A_j(U^\eps)\xi_j
=
\begin{pmatrix}
      v^\eps\cdot \xi & 0& \frac{1}{2}\RE a^\eps \,^{t}\xi \\ 
     0 &  v^\eps\cdot \xi & \frac{1}{2}\IM a^\eps\,^{t}\xi \\ 
     \ 2  \RE a^\eps\, \xi\ 
     &\ 2  \IM a^\eps\, \xi\  &\   v^\eps\cdot \xi {\rm I}_d 
\end{pmatrix},
\end{equation*}
and
\begin{equation*}
  L = \left(
    \begin{array}[l]{ccccc}
      
   0  &-\Delta &0& \dots & 0   \\
   \Delta  & 0 &0& \dots & 0  \\
   0& 0 &&0_{d\times d}& \\
   \end{array}
\right).
\end{equation*}
To be precise, the system is made symmetric thanks to the constant symmetrizer
\begin{equation*}
  S=\left(
    \begin{array}[l]{cc}
     {\rm I}_2 & 0\\
     0& \frac{1}{4}{\rm I}_d
    \end{array}
\right).
\end{equation*}
Once $v^\eps$ is known, one recovers $\phi^\eps$ by integrating in
time the first equation in \eqref{eq:Grenier},
\begin{equation*}
  \phi^\eps(t,x) = \phi_0(x)-\frac{1}{2}\int_0^t |v^\eps(s,x)|^2ds -
  \int_0^t|a^\eps(s,x)|^2ds, 
\end{equation*}
and since $\d_t (v^\eps-\nabla\phi^\eps)=0$, $v^\eps=\nabla
\phi^\eps$. 
Assuming that $a_0,\nabla \phi_0\in H^s(\R^d)$ for $s$ large (we will
always assume $a_0,\phi_0\in C_0^\infty(\R^d)$ in the forthcoming
applications), the limit $\eps\to 0$ leads to an asymptotic expansion
of the form
\begin{equation*}
  \phi^\eps \sim \phi+\eps \phi^{(1)}+\eps^2\phi^{(2)}+\dots,\quad
  a^\eps \sim a+\eps a^{(1)}+\eps^2 a^{(2)}+\dots
\end{equation*}
The leading order term is obtained by simply setting $\eps=0$ in
\eqref{eq:symL}: 
\begin{equation}
  \label{eq:euler-sym}
  \left\{
    \begin{aligned}
      &\d_t \phi  +\frac{1}{2}|\nabla
     \phi|^2+|a|^2=0,\quad \phi_{\mid t=0} =\phi_0,\\
&\d_t a+\nabla\phi\cdot \nabla
a+\frac{1}{2}a\Delta \phi =0,\quad a_{\mid t=0}=a_0.
    \end{aligned}
\right.
\end{equation}
Working with the intermediary unknown $v=\nabla \phi$, we get a system
of the form
\begin{equation*}
  \d_t U+\sum_{j=1}^dA_j(U)\d_j U=0,
\end{equation*}
and we infer the following result from
\cite{MUK86}:
\begin{proposition}\label{prop:MUK}
 Let $a_0,\phi_0\in C_0^\infty(\R^d)$,  with $\supp a_0,\supp
 \phi_0\subset K$. There exists
  $T_*>0$ and a unique solution $(\phi,a)\in C([0,T_*],H^\infty(\R^d))^2$ to
  \eqref{eq:euler-sym}. Moreover, $(\phi,a)$ remains compactly
  supported for $t\in [0,T_*]$, and 
  \begin{equation*}
       \supp \phi(t,\cdot),\ \supp a(t,\cdot)\subset K.
  \end{equation*}
\end{proposition}
The first part of the statement is a consequence of classical theory
for symmetric hyperbolic systems (see e.g. \cite{AlGe07,Majda}). The
property stated that initial compactly supported condition lead to a
zero speed of propagation is due to the structure of this hyperbolic
system, and is well understood from the simplest model of the Burgers
equation
\begin{equation*}
  \d_t u + u\d_xu=0,\quad u_{\mid t=0}=u_0\in C_0^\infty(\R).
\end{equation*}
Suppose we have a smooth solution on some time interval $[0,T_*]$. In particular,
\begin{equation*}
  \int_0^{T_*}\|\d_x u(t)\|_{L^\infty}dt<\infty. 
\end{equation*}
We have directly, for all $(t,x)\in [0,T_*]\times\R$,
\begin{equation*}
  |\d_tu(t,x)|\le \|\d_x u(t)\|_{L^\infty} |u(t,x)|.
\end{equation*}
Gronwall lemma then shows that if $u_0(x_0)=0$, then $u(t,x_0)=0$ for
all $t\in [0,T_*]$, hence the zero speed of propagation for smooth
solutions. As the matrix $A(U,\xi)$ is linear in $U$, the result
follows in the setting of \eqref{eq:euler-sym}. Note that to prove
this zero speed of propagation, we do not invoke the symmetry of $A$:
it was used in order to get Sobolev estimates (which ensure that $U\in
L^1([0,T_*],W^{1,\infty})$), but only the fact that it is (at least)
linear in $U$ is used at this stage.
We then have, for the same $T_*$ as in Proposition~\ref{prop:MUK}:
\begin{proposition}
Let $a_0,\phi_0\in C_0^\infty(\R^d)$.  There exists $T_*>0$ independent of $\eps\in ]0,1]$ such that for all
  $s\ge 0$, there exists $C=C(s)$ such that
  \begin{equation*}
    \|\phi^\eps-\phi\|_{L^\infty([0,T_*],H^s(\R^d))}+
    \|a^\eps-a\|_{L^\infty([0,T_*],H^s(\R^d))}\le C\eps. 
  \end{equation*}
\end{proposition}
To infer the pointwise description of $u^\eps$ at leading order, we
must in addition know $\phi^\eps$ up to $o(\eps)$, which is achieved
by considering the linearization of \eqref{eq:euler-sym}. At the next
step of the WKB expansion, we find that  
 \begin{equation*}
    \|\phi^\eps-\phi-\eps \phi^{(1)}\|_{L^\infty([0,T_*],H^s(\R^d))}+
    \|a^\eps-a-\eps a^{(1)}\|_{L^\infty([0,T_*],H^s(\R^d))}\le C\eps^2, 
  \end{equation*}
where the first corrector $(\phi^{(1)},a^{(1)})$ solves the system:
\begin{equation*}
  \left\{
    \begin{aligned}
     &\d_t \phi^{(1)} +\nabla \phi\cdot \nabla \phi^{(1)} + 2
\RE\( \overline a a^{(1)}\)=0,\\
& \d_t a^{(1)}+\nabla \phi\cdot \nabla a^{(1)} + \nabla \phi^{(1)}\cdot
\nabla a +\frac{1}{2} a^{(1)}\Delta \phi + \frac{1}{2} a\Delta
\phi^{(1)}= \frac{i}{2}\Delta a, \\
 & \phi^{(1)}_{\mid t=0}=0\quad ;\quad a^{(1)}_{\mid t=0}=0.
    \end{aligned}
\right.
\end{equation*}
At higher order $k\ge 2$, the corrector $(\phi^{(k)},a^{(k)})$ is given by:
\begin{equation*}
  \left\{
    \begin{aligned}
     &\d_t \phi^{(k)} +\nabla \phi\cdot \nabla \phi^{(k)} + 2
\RE\( \overline a a^{(k)}\)=F_k\(\(\nabla \phi^{(\ell)}\)_{1\le \ell\le
  k-1},\(a^{(\ell)}\)_{1\le \ell\le
  k-1} \) ,\\
& \d_t a^{(k)}+\nabla \phi\cdot \nabla a^{(k)} + \nabla \phi^{(k)}\cdot
\nabla a +\frac{1}{2} a^{(k)}\Delta \phi + \frac{1}{2} a\Delta
\phi^{(k)}= \frac{i}{2}\Delta a^{(k-1)}\\
&\phantom{\d_t a^{(j)}+\nabla \phi\cdot \nabla a^{(j)} } -\sum_{1\le \ell\le k-1}\nabla \phi^{(\ell)}\cdot \nabla a^{(k-\ell)}-
\frac{1}{2}\sum_{1\le \ell\le k-1}a^{(k-\ell)}\Delta \phi^{(\ell)},\\
 & \phi^{(k)}_{\mid t=0}=0\quad ;\quad a^{(k)}_{\mid t=0}=0,
    \end{aligned}
\right.
\end{equation*}
for some function $F_k$ which is a polynomial in its arguments,
without constant term, and whose
precise expression is unimportant here. The left hand side is always
the linearization of the left hand side of \eqref{eq:euler-sym} about
$(\phi,a)$, and the right hand side depends on previous correctors. We
infer (see \cite{Grenier98,CaBook2}), by induction:
\begin{proposition}\label{prop:BKW-gen1}
  Let $a_0,\phi_0\in C_0^\infty(\R^d)$.  Let $T_*>0$ given by
  Proposition~\ref{prop:MUK}. For all $k\ge 1$, there exists a unique
  solution $(\phi^{(k)},a^{(k)})\in C([0,T_*],H^\infty(\R^d))^2$ to the
  above system, and for all
  $s\ge 0$, there exists $C=C(k,s)$ such that
  \begin{align*}
    &\left\|\phi^\eps-\phi-
      \eps\phi^{(1)}-\ldots-\eps^k\phi^{(k)}\right\|_{L^\infty([0,T_*],H^s(\R^d))}\\
&+ 
    \left\|a^\eps-a-\eps
      a^{(1)}-\ldots-\eps^ka^{(k)}\right\|_{L^\infty([0,T_*],H^s(\R^d))}\le
    C\eps^{k+1}.  
  \end{align*}
In addition, if $\supp a_0,\supp
 \phi_0\subset K$, then
 $(\phi^{(k)},a^{(k)})$ remains compactly 
  supported for $t\in [0,T_*]$, and 
  \begin{equation*}
       \supp \phi^{(k)}(t,\cdot),\ \supp a^{(k)}(t,\cdot)\subset K.
  \end{equation*}
\end{proposition}
The support property is a consequence of the same argument as in the
proof of Proposition~\ref{prop:MUK}. Using the embedding
$H^s(\R^d)\subset L^\infty(\R^d)$ for $s>d/2$, we also deduce from the above
error estimate the bound, for $k\ge 1$:
\begin{equation}\label{eq:error-j}
  \left\| u^\eps - \(\sum_{\ell=0}^{k-1} \eps^\ell
    a^{(\ell)}\)\exp\(\frac{i}{\eps}\sum_{\ell=0}^k \eps^\ell \phi^{(\ell)}
    \)\right\|_{L^\infty([0,T_*],L^2\cap L^\infty(\R^d))} =\O\(\eps^{k}\),
\end{equation}
with the convention $(\phi^{(0)},a^{(0})=(\phi,a)$. The standard form
of WKB expansions,
\begin{equation*}
  u^\eps(t,x) =\( {\tt a}(t,x)+\eps {\tt a}_1(t,x)+\ldots +\eps^k{\tt
    a}_k(t,x)\)e^{i\phi(t,x)/\eps}+
\O_{L^\infty_{T_*}(L^2\cap L^\infty)}\(\eps^{k+1}\),
\end{equation*}
is then obtained by setting
\begin{equation*}
  {\tt a}= a e^{i\phi^{(1)}},\quad {\tt a}_1 = a^{(1)}
  e^{i\phi^{(1)}}+ia \phi^{(2)}e^{i\phi^{(1)}},\quad \text{etc.}
\end{equation*}
\begin{remark}[Higher order nonlinearities]
  If instead of \eqref{eq:mono}, one considers
  \begin{equation*}
 i\eps\d_t u^\eps+\frac{\eps^2}{2}\Delta u^\eps=
  |u^\eps|^{2\si}u^\eps\quad ;\quad u^\eps_{\mid t=0}=a_0e^{i\phi_0/\eps},
\end{equation*}
with $\sigma\ge 2$ an integer, then the justification of WKB analysis
requires a different approach. We refer to \cite{ACARMA,ChironRousset}
for two different 
proofs, which show that the conclusions of the propositions stated in
this section remain valid. 
\end{remark}

\begin{remark}[Focusing nonlinearity]
  If instead of \eqref{eq:mono}, one considers a cubic focusing nonlinearity,
  \begin{equation*}
 i\eps\d_t u^\eps+\frac{\eps^2}{2}\Delta u^\eps=
 - |u^\eps|^{2}u^\eps\quad ;\quad u^\eps_{\mid t=0}=a_0e^{i\phi_0/\eps},
\end{equation*}
then the analogue of \eqref{eq:euler-sym} is no longer hyperbolic, but
elliptic. Working with analytic initial data $(\phi_0,a_0)$ is then
necessary in order to solve \eqref{eq:euler-sym}
(\cite{LeNgTe18,GuyCauchy}), and this 
is a framework where nonlinear WKB analysis is fully justified
(\cite{PGX93,ThomannAnalytic}). However, analyticity is incompatible with an initial
compact support. On the other hand, in the weakly nonlinear case
$\gamma =1$ (and more generally if $\gamma\ge 1$), it is possible to
justify WKB analysis with a focusing nonlinearity and compactly
supported initial data (see e.g. \cite{CaBKW} or \cite[Chapter~2]{CaBook2}). 
\end{remark}

\section{The monokinetic case with a potential}
\label{sec:monopot}

In this section, we first recall some elements of WKB analysis in the
linear case. We then show how this case can be merged with the analysis
presented in the previous section, when $\gamma=0$. We sketch how the
case of a weaker nonlinearity, $0<\gamma<1$. To conclude, we briefly discuss the weakly nonlinear
r\'egime $\gamma=1$, and more generally the situation $\gamma\ge 1$. 

 \subsection{Linear case}
\label{sec:linear}

 The eikonal
 equation associated to
 \begin{equation}\label{eq:Slin}
   i\eps\d_t u^\eps+\frac{\eps^2}{2}\Delta u^\eps = V u^\eps\quad
   ;\quad u^\eps_{\mid t=0} = a_0e^{i\varphi_0/\eps},
 \end{equation}
 that is, without initial rapid oscillation, reads:
\begin{equation}
  \label{eq:eik}
  \d_t \phi_{\rm eik}+\frac{1}{2}|\nabla \phi_{\rm eik}|^2+V=0\quad
  ;\quad \phi_{{\rm eik}\mid t=0}=\varphi_0.
\end{equation}
In this subsection, we assume that $\varphi_0$ is smooth and at most
quadratic, in the same sense as in \eqref{eq:V}. 
This eikonal equation is solved by introducing the classical
trajectories, solving
\begin{equation}\label{eq:Hamilton}
  \dot x(t,y) = \xi(t,y),\quad x(0,y)=y\quad ;\quad \dot \xi(t,y) =
  -\nabla V\(x(t,y)\),\quad \xi(0,y)=\nabla \varphi_0(y).
\end{equation}
As $V$ is at most quadratic, from \eqref{eq:V}, the above system has a
unique, global, smooth solution, and in addition
\begin{equation*}
  \nabla_y x(t,y) = {\rm I}_d +\O(t), 
\end{equation*}
uniformly in $y\in \R^d$, for any matricial norm on $\R^{d\times
  d}$. Therefore, the Jacobi determinant
\begin{equation*}
  J_t(y) = \operatorname{det} \nabla_y x(t,y) ,
\end{equation*}
remains non-zero and bounded on some time interval $[0,T]$ with
$T>0$. Since we also have, by uniqueness in ordinary differential
equations,
\begin{equation*}
  \nabla \phi_{\rm eik} \(t,x(t,y)\)= \xi(t,y),
\end{equation*}
for any smooth solutions to \eqref{eq:eik}, the global inversion
theorem implies the following result (see also \cite[Proposition~1.9]{CaBook2}):
\begin{lemma}\label{lem:eik}
Let  $V$ satisfying \eqref{eq:V}, and $\varphi_0$ satisfying the same
properties.  There exists $T>0$ and a 
  unique solution $\eik\in C^\infty \([0,T]\times \R^d\)$ to
  \eqref{eq:eik}. 
In addition, this solution is at most quadratic in space: $\partial_x^\alpha
  \eik \in L^\infty([0,T]\times\R^d)$ as soon as $|\alpha|\ge
  2$. There exists $C>1$ such that the Jacobi determinant satisfies:
  \begin{equation*}
    \frac{1}{C}\le J_t(y)\le C,\quad \forall (t,y)\in [0,T]\times
    \R^d. 
  \end{equation*}
\end{lemma}
Since the above relations imply
\begin{equation*}
  \dot x(t,y) = \nabla\eik \(t,x(t,y)\),
\end{equation*}
we infer the classical formula
\begin{equation}\label{eq:dtJ}
  \d_t J_t(y) = J_t(y)\Delta \eik\(t,x(t,y)\).
\end{equation}
In the linear case, the leading order amplitude is given by the linear
transport equation
\begin{equation*}
  \d_t a + \nabla\eik\cdot\nabla a
  +\frac{1}{2}a\Delta\eik=0\quad;\quad a_{\mid t=0}=a_0.
\end{equation*}
Following the classical trajectories, this transport equation becomes
trivial, since $A(t,y):= \sqrt{J_t(y)}a\(t,x(t,y)\)$ satisfies $\d_t
A=0$. 
\subsection{Supercritical case: $\gamma=0$}

We consider the same framework as in the previous section, now with a potential:
\begin{equation}
  \label{eq:monopot}
  i\eps\d_t u^\eps+\frac{\eps^2}{2}\Delta u^\eps=
  Vu^\eps+|u^\eps|^2u^\eps\quad ;\quad u^\eps_{\mid t=0}=a_0e^{i\phi_0/\eps}.
\end{equation}
 As noticed in \cite{CaBKW}, it is possible to adapt the above WKB
 analysis in the presence of an external potential satisfying
 \eqref{eq:V} by simply mixing the standard approach followed in the linear
 case (see e.g. \cite{Robert}) and Grenier's method.

\subsubsection{Introducing the nonlinearity}

As noticed in \cite{CaBKW}, the approach presented in the case $V=0$
for the nonlinear case can be adapted by changing the representation
\eqref{eq:reprGrenier} to
\begin{equation*}
  u^\eps(t,x) = a^\eps(t,x)e^{i\eik(t,x)/\eps+i\phi^\eps(t,x)/\eps},
\end{equation*}
where $\eik$ solves \eqref{eq:eik} with $\varphi_0\equiv 0$, and requiring
\begin{equation}\label{eq:systcmp2}
  \left\{
    \begin{aligned}
      & \d_t \phi^\eps +\nabla \eik\cdot\nabla \phi^\eps+
      \frac{1}{2}|\nabla \phi^\eps|^2+
      |a^\eps|^2 = 0,\\
      & \d_t a^\eps +\nabla \eik\cdot\nabla a^\eps+ \nabla \phi^\eps
      \cdot \nabla a^\eps 
    +\frac{1}{2}a^\eps \Delta \eik  +\frac{1}{2}a^\eps \Delta
      \phi^\eps =i\frac{\eps}{2}\Delta a^\eps ,\\
& \quad \phi^\eps_{\mid t=0}=\phi_0\quad ;\quad   a^\eps_{\mid
      t=0}=a_0. 
    \end{aligned}
\right.
\end{equation} 
The new terms compared to \eqref{eq:Grenier} involve $\eik$, and since
$\eik$ is at most quadratic in space, it turns out that they can be
estimated like (semilinear) perturbative terms (using commutator
estimates for the transport part). The natural limit for
\eqref{eq:systcmp2} when $\eps\to 0$ is given by
\begin{equation}\label{eq:systcmp-lim}
  \left\{
    \begin{aligned}
      & \d_t \phi +\nabla \eik\cdot\nabla \phi+
      \frac{1}{2}|\nabla \phi|^2+
      |a|^2 = 0,\\
      & \d_t a +\nabla \eik\cdot\nabla a+ \nabla \phi
      \cdot \nabla a
    +\frac{1}{2}a \Delta \eik  +\frac{1}{2}a \Delta
      \phi =0 ,\\
& \quad \phi_{\mid t=0}=\phi_0\quad ;\quad   a_{\mid
      t=0}=a_0. 
    \end{aligned}
\right.
\end{equation} 
The following result is a consequence of \cite{CaBKW}:
\begin{proposition}\label{prop:BKWNLpot}
  Let $a_0,\phi_0\in C_0^\infty(\R^d)$, $V$ satisfying \eqref{eq:V},
  and $T$, $\eik$ given by Lemma~\ref{lem:eik}.   There exists $0<T_*\le T$
  independent of $\eps\in ]0,1]$ such that \eqref{eq:systcmp2} has a
  unique solution $(\phi^\eps,a^\eps)\in C([0,T_*],H^\infty(\R^d))^2$,
  \eqref{eq:systcmp-lim} has a 
  unique solution $(\phi,a)\in C([0,T_*],H^\infty(\R^d))^2$, and 
for all 
  $s\ge 0$, there exists $C=C(s)$ such that 
  \begin{equation*}
    \|\phi^\eps-\phi\|_{L^\infty([0,T],H^s(\R^d))}+
    \|a^\eps-a\|_{L^\infty([0,T],H^s(\R^d))}\le C\eps. 
  \end{equation*}
\end{proposition}
The correctors $\(\phi^{(j)},a^{(j)}\)_{j\ge 1}$ as obtained in the
same fashion as in Section~\ref{sec:mono}. The only difference is that
the operator $\d_t$ is replaced by
\begin{equation*}
  \d_t + \nabla\eik\cdot\nabla +\frac{1}{2}\Delta \eik.
\end{equation*}

\subsubsection{Finite speed of propagation: following the classical
  trajectories}

In order to prove that if $a_0\in C_0^\infty(\R^d)$, the solution to
\eqref{eq:Slin} remains compactly supported in the support of $a_0$
transported by the classical flow \eqref{eq:Hamilton}, it is standard
to introduce the following change of unknown function
(e.g. \cite{Robert,CaBook2}):
\begin{equation*}
  A(t,y) := \sqrt{J_t(y)} a\(t,x(t,y)\),
\end{equation*}
where $a$ solves the transport equation
\begin{equation*}
  \d_t a+ \nabla\eik\cdot\nabla a+\frac{1}{2}a\Delta \eik=0 \quad;\quad
  a_{\mid t=0}=a_0.
\end{equation*}
as given by WKB analysis. Indeed, using \eqref{eq:dtJ}, we easily
check that $A$ is constant in time, $\d_t A=0$. Correctors
$(a^{(k)})_{k\ge 1}$ in the
(linear) WKB analysis solve the equation
\begin{equation*}
  \d_t a^{(k)}+ \nabla\eik\cdot\nabla a^{(k)}+\frac{1}{2}a^{(k)}\Delta
  \eik=\frac{i}{2}\Delta a^{(k-1)} \quad;\quad
  a^{(k)}_{\mid t=0}=0,
\end{equation*}
with the convention $a^{(0)}=a$. Setting
\begin{equation*}
    A^{(k)}(t,y) := \sqrt{J_t(y)} a^{(k)}\(t,x(t,y)\),
  \end{equation*}
  we infer that
  \begin{equation*}
    \supp A^{(k)}(t,\cdot)\subset \supp a_0,\quad \forall t\in [0,T],\
    \forall k\ge 0,
  \end{equation*}
  where $T$ is given by Lemma~\ref{lem:eik}. Thus, for $t\in [0,T]$, up to
  $\O(\eps^\infty)$, $u^\eps$ remains compactly supported, in the
  support of $a_0$ transported by the classical flow.
  \smallbreak

  In the nonlinear case, we check that the same argument remains
  valid. Consider $\eik$ solution to \eqref{eq:eik}, and $(\phi,a)$
  solving \eqref{eq:systcmp-lim}. The natural adaptation of the above
  computation consists in showing that if $\phi_0,a_0\in
  C_0^\infty(\R^d)$, the new unknown $(\psi,A)$,
  defined by
  \begin{equation}\label{eq:redress}
      A(t,y) := \sqrt{J_t(y)} a\(t,x(t,y)\),\quad \psi(t,y) := \phi \(t,x(t,y)\),
  \end{equation}
  enjoys a zero speed of propagation. Note that in view of
  Proposition~\ref{prop:BKWNLpot}, we already know that $\phi,a\in
  C([0,T_*],H^\infty(\R^d))$, so it suffices to check that $(\psi,A)$
  solves a system for which the argument presented on the toy model of
  Burgers equation in Section~\ref{sec:mono} remains valid.
  Introducing
  \begin{equation*}
    M(t,y) = \nabla_y x(t,y)\in \R^{d\times d},
  \end{equation*}
whose determinant is by definition $J_t(y)$, we find:
\begin{align*}
  &\d_t \psi +\frac{1}{2}
   \<M^{-1}\nabla \psi, M^{-1}\nabla \psi\>+
    \frac{1}{J_t(y)}|A|^2 =0 ,\quad \psi_{\mid t=0}=\phi_0,\\
  &\d_t A= -\sqrt{J_t(y)}\( \nabla \phi\cdot \nabla
    a+\frac{1}{2}a\Delta\phi\)\(t,x(t,y)\),\quad A_{\mid t=0}=a_0.
\end{align*}
We do not express the right hand side of the last equation in terms of
$(\psi,A)$: differentiating the first equation with respect to $y$,
the bounds stated in Proposition~\ref{prop:BKWNLpot} make it possible
to infer an inequality of the form
\begin{equation*}
  |\d_t\nabla \psi(t,y)|+|\d_tA(t,y)|\lesssim |\nabla
  \psi(t,y)|+|A(t,y)|,\quad (t,y)\in [0,T_*]\times \R^d.
\end{equation*}
Therefore, if $\supp \phi_0,\supp a_0\subset K$, then
$\supp\nabla \psi(t,\cdot),\supp A(t,\cdot)\subset K$ for all $t\in
[0,T_*]$. Integrating in time the equation solved by $\psi$, we
conclude to the zero speed of propagation for $(\psi,A)$. Arguing like
in Section~\ref{sec:mono} for the correctors, we have:
\begin{proposition}\label{prop:zerospeedpot}
  Let $\phi_0,a_0\in C_0^\infty(\R^d)$ with  $\supp \phi_0,\supp
  a_0\subset K$. There for any $t\in [0,T_*]$, where $T_*$ is given by
  Proposition~\ref{prop:BKWNLpot},
  \begin{equation*}
    \supp\psi(t,\cdot),\supp A(t,\cdot)\subset K,
  \end{equation*}
  where $\psi$ and $A$ are related to $\phi$ and $a$ through
  \eqref{eq:redress}. The same is true for the correctors
  $(\psi^{(k)},A^{(k)})_{k\ge 1}$ corresponding to the next terms
  $(\phi^{(k)},a^{(k)})_{k\ge 1}$  in the
  asymptotic expansion in \eqref{eq:systcmp2}. 
\end{proposition}

  \begin{remark}[Special potentials]
    In the case where $V$ is linear in $x$ or isotropic quadratic,
    explicit formulas allow to bypass the above arguments. If
    $V(x)=E\cdot x$ for some (constant) $E\in \R^d$ and $u^\eps$
    solves \eqref{eq:NLS}, then
    \begin{equation*}
      v^\eps(t,x) = u^\eps\(t,x-\frac{t^2}{2}E\)e^{i\(tE\cdot x
        -\frac{t^3}{3}|E|^2\)/\eps} 
    \end{equation*}
    solves \eqref{eq:mono}.  If $V(x) =\frac{\omega^2}{2}|x|^2$,
    $\omega>0$, then
    \begin{equation*}
      w^\eps(t,x) = \frac{1}{\(1+(\omega t)^2\)^{d/4}}
      u^\eps\(\frac{\operatorname{arctan} (\omega
        t)}{\omega},\frac{x}{\sqrt{1+(\omega   t)^2}}\)
      e^{i\frac{\omega^2t}{1+(\omega t)^2}\frac{|x|^2}{2\eps}}
      \end{equation*}
      solves
      \begin{equation*}
        i\eps\d_t w^\eps + \frac{\eps^2}{2}\Delta w^\eps =
        (1+t^2)^{d/2-1}|w^\eps|^2w^\eps\quad ;\quad w^\eps_{\mid t=0}
        = a_0e^{i\phi_0/\eps}. 
      \end{equation*}
      If $d=2$ (the cubic nonlinearity is $L^2$-critical), we recover
      exactly \eqref{eq:mono}. Otherwise, a (smooth) time dependent factor has
      appeared, which obviously does not change the conclusion of
      Propositions~\ref{prop:MUK} and \ref{prop:BKW-gen1}. The case of
      a potential with the opposite sign is obtained by changing
      $\omega$ to $i\omega$ in the formulas. See
      e.g. \cite[Section~11.2]{CaBook2} and references therein regarding these
      changes of unknown functions. For such potentials, the classical
      trajectories given by \eqref{eq:Hamilton} are computed
      explicitly, and we can check directly the conclusions of
      Proposition~\ref{prop:zerospeedpot}. 
  \end{remark}

\subsection{Weaker nonlinearity}

We now consider the case $0<\gamma<1$. This case is still a
supercritical case as far as WKB analysis is concerned, in the sense
described in the introduction: a ``natural'' asymptotic expansion of
the solution $u^\eps$ still involves a system of equations which is
not closed. As noticed in \cite{CaBKW}, this intermediary case can be
handled like the case $\gamma=0$, by replacing \eqref{eq:systcmp2}
with 
\begin{equation}\label{eq:systcmp2}
  \left\{
    \begin{aligned}
      & \d_t \phi^\eps +\nabla \eik\cdot\nabla \phi^\eps+
      \frac{1}{2}|\nabla \phi^\eps|^2+
      \eps^\gamma|a^\eps|^2 = 0,\\
      & \d_t a^\eps +\nabla \eik\cdot\nabla a^\eps+ \nabla \phi^\eps
      \cdot \nabla a^\eps 
    +\frac{1}{2}a^\eps \Delta \eik  +\frac{1}{2}a^\eps \Delta
      \phi^\eps =i\frac{\eps}{2}\Delta a^\eps ,\\
& \quad \phi^\eps_{\mid t=0}=\phi_0\quad ;\quad   a^\eps_{\mid
      t=0}=a_0. 
    \end{aligned}
\right.
\end{equation} 
The matrices $A_j$ and $S$ now depend on $\eps$, in an explicit way,
and the asymptotic expansion of $(\phi^\eps,a^\eps)$ involves more
terms. Let $N= [1/\gamma]$, where $[r]$ is the largest integer not larger than
$r>0$: $N$ new intermediary terms appear compared to the case
$\gamma=0$,
\begin{align*}
  &\phi^\eps = \eps^\gamma \tilde \phi + \eps^{2\gamma} \tilde
  \phi^{(1)}+\dots +\eps^{N\gamma}\tilde \phi^{(N)}+\O(\eps),\\
&a^\eps = a+ \eps^\gamma
  a^{(1)}+\dots +\eps^{N\gamma}a^{(N)}+\O(\eps),
\end{align*}
where the estimate holds in $L^\infty([0,T],H^s)$ for any
$s>0$. This can be seen by  setting $\tilde \phi^\eps =
\eps^{-\gamma}\phi^\eps$: the leading order term is given by 
\begin{equation*}
\left\{
\begin{aligned}
    \partial_t \widetilde \phi  +\nabla 
    \phi_{\rm eik}\cdot \nabla \widetilde \phi+ 
    |a|^2&= 0\ ; \quad \widetilde \phi_{\mid t=0}=0, \\
    \partial_t a +\nabla \phi_{\rm eik} \cdot \nabla a +\frac{1}{2}a 
\Delta \phi_{\rm eik} &= 0\ ; \quad a_{\mid t=0}= a_0 .
\end{aligned} 
\right.
\end{equation*}
The leading order amplitude solves the same transport equation as in
the linear case, and it is readily observed that the analogue of
Proposition~\ref{prop:zerospeedpot} remains valid, up to adapting the
hierarchy of equations. 
\subsection{Weakly nonlinear and linearizable cases}

We now assume $\gamma\ge 1$. As in \cite{CaBKW} (or
\cite[Chapter~2]{CaBook2}), we present a strategy for any $\gamma\ge
1$, and emphasize the fact that the value $\gamma=1$ is specific. In
this setting, the coupling between phase and amplitude changes
dramatically: rapid oscillations are described by $\phi_{\rm eik}$
only, and the analysis consists in expanding the amplitude $a^\eps =
u^\eps e^{-i\eik/\eps}$ in powers of $\eps$:
\begin{equation*}
  \d_t a^\eps +\nabla \eik\cdot\nabla a^\eps
    +\frac{1}{2}a^\eps \Delta \eik  =i\frac{\eps}{2}\Delta a^\eps -i\eps^{\gamma-1}|a^\eps|^2a^\eps\quad;
 \quad  a^\eps_{\mid
      t=0}=a_0. 
\end{equation*}
Like above, the case when $\gamma>1$ is not an integer requires a
special asymptotic expansion, and we do not discuss this case. When
$\gamma>1$, the leading order amplitude satisfies the same transport
equation as in the linear case. When $\gamma=1$, it satisfies
\begin{equation*}
  \d_t a+\nabla \eik\cdot\nabla a
    +\frac{1}{2}a \Delta \eik  =-i|a|^2a\quad;
 \quad  a_{\mid
      t=0}=a_0. 
\end{equation*}
Following the classical trajectories, that is resuming the change of
unknown function \eqref{eq:redress}, this equation reads
$\d_tA =-i J_t(y)^{-1}|A|^2A$. In particular, $\d_t|A|^2=0$, and the
nonlinear effect in $a$ consists of a phase selfmodulation. In
particular, the support of $A(t,\cdot)$ is independent of $t\in
[0,T]$. The same is true for all correctors in the asymptotic
expansion, as can be checked easily.

\section{Separation of states}
\label{sec:separation}

We complete the proof of Theorem~\ref{theo:superp}, by proving the
nonlinear superposition. 
For $1\le j\le N$, let $\chi_j\in C_0^\infty(\R^d,\R)$, $0\le \chi_j\le 1$, with
\begin{equation*}
  \chi_j\equiv  1\text{ on }\supp \alpha_j,\quad \supp \chi_{j_1}\cap \supp
  \chi_{j_2}=\emptyset\text{ if }j_1\not =j_2. 
\end{equation*}
We set
\begin{equation*}
  a_0(x) =\sum_{j=1}^N \alpha_j(x),\quad \phi_0(x) =\sum_{j=1}^N 
  \varphi_j(x)\chi_j(x). 
\end{equation*}
Then $a_0,\phi_0\in C^\infty_0(\R^d)$, $\phi_0$ is real-valued, and
\begin{equation*}
  u_0^\eps(x) = a_0(x)e^{i\phi_0(x)/\eps}. 
\end{equation*}
We can then resume the analysis from the monokinetic case as presented
in Sections~\ref{sec:mono} and \ref{sec:monopot}, with the same
notations. Let $\eik$ be given by Lemma~\ref{lem:eik} (it does not
depend on the initial data, but only on $V$).

\subsection{Supercritical case}

 When $\gamma=0$, the WKB analysis for
each $u_j^\eps$, solution to 
\eqref{eq:NLSj}, involves the following system:
\begin{equation}
  \label{eq:euler-sym-j}
  \left\{
    \begin{aligned}
      &\d_t \phi_j  +\nabla \eik\cdot \nabla \phi_j+\frac{1}{2}|\nabla
     \phi_j|^2+|a_j|^2=0,\quad \phi_{j\mid t=0} =\varphi_j\chi_j,\\
&\d_t a_j+\nabla\eik\cdot
\nabla a+\frac{1}{2}a\Delta\eik+\nabla\phi_j\cdot \nabla 
a_j+\frac{1}{2}a_j\Delta \phi_j =0,\quad a_{j\mid t=0}=\alpha_j.
    \end{aligned}
\right.
\end{equation}
To simplify the discussion, suppose first
that $V=0$, hence $\eik =0$. Each solution to \eqref{eq:euler-sym-j} remains smooth on some
time interval $[0,T_j]$ for some $0<T_j\le T$, and, on this time interval,
enjoys a zero speed of propagation. 
As a consequence of Proposition~\ref{prop:MUK}, we have
\begin{equation*}
  \phi = \sum_{j=1}^N\phi_j,\quad a=\sum_{\j=1}^N a_j,
\end{equation*}
since nonlinear terms containing two indices $j_1\not =j_2$ involve two
functions whose supports are disjoint. Also, for all $k\ge 1$, the
correctors satisfy
\begin{equation*}
  \phi^{(k)} = \sum_{j=1}^N\phi_j^{(k)},\quad a=\sum_{\j=1}^N a_j^{(k)}.
\end{equation*}
Set
\begin{equation*}
  T^*=\min \(T_*,T_1,\ldots,T_N\).
\end{equation*}
As we have, in view of
\eqref{eq:error-j}, in $L^\infty([0,T^*],L^2\cap L^\infty)$, for any $k>0$, 
\begin{align*}
  &u^\eps - \(\sum_{\ell=0}^{k-1} \eps^\ell
    a^{(\ell)}\)\exp\(\frac{i}{\eps}\sum_{\ell=0}^k \eps^\ell
    \phi^{(\ell)}\)=\O\(\eps^k\),\\
&u_j^\eps - \(\sum_{\ell=0}^{k-1} \eps^\ell
    a_j^{(\ell)}\)\exp\(\frac{i}{\eps}\sum_{\ell=0}^k \eps^\ell
    \phi_j^{(\ell)}\)=\O\(\eps^k\),\quad j=1,\dots,N,
\end{align*}
we obtain Theorem~\ref{theo:superp} in the case $V=0$.
In the case where $V$ is not trivial, we just have to resume the above
arguments by replacing the functions $(\phi,a)$ (possibly with indices
and/or superscripts) with $(\psi,A)$, as defined by the change of
unknown function \eqref{eq:redress}, which involves only $V$ (see
\eqref{eq:Hamilton}), and not
the initial data. 

\subsection{Other cases}
When $\gamma>0$, we have seen that the leading order amplitude is the same
as in the linear case, up to a phase modulation.  Leading
order oscillations are given by $\eik$, where we now set
$\varphi_0=\phi_0$ in \eqref{eq:eik}. The features used in
the supercritical case then remain, regarding the evolution of the
support of the terms involved in WKB analysis.

\section{On the role of the cutoff function(s)}
\label{sec:chi}

\subsection{WKB analysis for the linear Schr\"odinger equation}

Consider \eqref{eq:Slin}  in the presence of rapid initial
oscillations,

\begin{equation}
  \label{eq:Slin-phase}
  i\eps \d_t u^\eps+\frac{\eps^2}{2}\Delta u^\eps = V u^\eps\quad
  u^\eps_{\mid t=0} = a_0 e^{i\phi_0/\eps},
\end{equation}
with $a_0\in C_0^\infty(\R^d)$ and $\phi_0\in C^\infty(\R^d,\R)$. Like
in Section~\ref{sec:separation}, consider $\chi\in
C_0^\infty(\R^d,[0,1])$, with
\begin{equation*}
  \chi\equiv 1\quad\text{on }\supp a_0. 
\end{equation*}
For any such function $\chi$, we have $u^\eps_{\mid t=0} = a_0
e^{i\chi\phi_0/\eps}$. However, the eikonal equation now depends on
$\chi$, as \eqref{eq:eik} becomes
\begin{equation*}
  \d_t \eik + \frac{1}{2}|\nabla \eik|^2+V =0\quad ;\quad \phi_{{\rm
      eik}\mid t=0} = \chi \phi_0.
\end{equation*}
As recalled in Section~\ref{sec:linear} (in the case $\phi_0=0$), the
solution is constructed, locally in time, \emph{via} the classical
trajectories, or, equivalently, through characteristic curves. As $V$
is smooth, the slope of characteristic curves at time $t=0$ 
is uniformly bounded on the support of $a_0$. By finite speed of
propagation, there exists $T(\chi)>0$ such that $\eik$ does not depend
on $\chi$ for $t\in [0,T(\chi)]$. In practice, the introduction of
$\chi$ may shorten the time interval of validity of WKB analysis, as
we now illustrate.
\smallbreak

Let $d=1$, $V=0$, and $\phi_0(x)=x^2/2$. The solution to the eikonal
equation (without cutoff $\chi$) is given explicitly by
\begin{equation*}
  \eik(t,x) = \frac{x^2}{2(1+t)}. 
\end{equation*}
This is a case where there is no singularity for $t\ge 0$ (but a
caustic reduced to one point at $t=-1$). Indeed, the classical
trajectories, solving
\begin{equation*}
  \dot x(t,y) = \xi(t,y),\quad x(0,y)=y\quad ;\quad \dot \xi(t,y) =
  0,\quad \xi(0,y) = \phi_0'(y)=y,
\end{equation*}
are given by
\begin{equation*}
  x(t,y) = (1+t)y, 
\end{equation*}
obviously inverted, for all $t\ge 0$,  as
\begin{equation*}
  y(t,x)=\frac{x}{1+t},
\end{equation*}
and the leading order amplitude in WKB analysis is given by
\begin{equation*}
  a(t,x) = \frac{1}{\sqrt{1+t}}a_0\(\frac{x}{1+t}\). 
\end{equation*}
For $\chi$ a (usual) cutoff function as above, $\chi \phi_0$ has two
humps: in the presence of $\chi$, $y\mapsto x(t,y)$ ceases to be
invertible for all $t\ge 0$ ($\eik '$ solves the Burgers equation),
but for short time (independent of $\eps$, but depending on $\chi$), $a(t)e^{i\eik(t)/\eps}$
does not depend on $\chi$.

\subsection{Supercritical WKB analysis for the nonlinear Schr\"odinger
  equation}

In the case addressed in Section~\ref{sec:mono}, the above eikonal
equation is replaced by \eqref{eq:euler-sym}. By considering the
gradient of the phase instead of the phase,  the Burgers equation (in
the case of WKB analysis for the linear Schr\"odinger equation without
potential) is replaced by the symmetrization of the Euler
equation. Like above, finite speed of propagation implies that the
introduction of a cutoff function in the initial phase does not alter
the solution to \eqref{eq:euler-sym} on some time interval
$[0,T(\chi)]$, for some $T(\chi)>0$ possibly depending on $\chi$. This time
is of course independent of $\eps$, as $\eps$ is absent from
\eqref{eq:euler-sym}. This is why in Remark~\ref{rem:wigner1}, the
Wigner measure does not depend on the $\chi_j$'s, even though its
construction seems to depend on these cutoff functions: the time of
validity that we can prove may, on the other hand, depend on the
choice of these cutoff functions.
\smallbreak

We conclude this discussion by an illustration similar to the one
given in the previous subsection. Let $a_0\in C^\infty_0(\R^d)$, and
assume that for $s > d/2 +1$, $\left\|a_0\right\|_{H^s(\R^d)}$ is
sufficiently small. Suppose also that $v_0=\nabla \phi_0$ satisfies:
$\nabla^{2}v_0\in H^{s-1}(\R^d)$, $\nabla v_0\in
  L^\infty(\R^d)$, and there exists $\delta>0$ such that for all
  $x\in \R^d$, ${\rm dist}(\operatorname{Sp}(\nabla v_0(x)),\R_-)\ge \delta$, 
where we denote by $\operatorname{Sp}(M)$ the
spectrum of a matrix $M$. Then it follows from the main result in
\cite{Gra98} that \eqref{eq:euler-sym} has a global (in the future)
solution
\begin{equation*}
  a,v-\overline v \in C^j([0,\infty),H^{s-j}(\R^d)), \quad j=0,1,
\end{equation*}
where $\overline v $ is the unique, global smooth solution to the
(multidimensional) Burgers equation
\begin{equation*}
  \d_t \overline v + \overline v\cdot \nabla \overline v =0,\quad
  \overline v_{\mid t=0} = \nabla \phi_0. 
\end{equation*}
We may for instance consider $\phi_0(x)=|x|^2/2$ (see the previous
subsection), and then
\begin{equation*}
  \overline v(t,x) = \frac{x}{1+t}. 
\end{equation*}
On the other hand, if $\phi_0$ is multiplied by a cutoff function
$\chi$, then the initial data in \eqref{eq:euler-sym} belong to
$C_0^\infty(\R^d)$: it follows from \cite{MUK86} that the
corresponding solution develops a singularity in finite time. Like in
the previous subsection, the introduction of the cutoff $\chi$ reduces
the lifespan of the solution involved in WKB analysis but, for short
time, does not alter the asymptotic description of the solution
$u^\eps$.

\section{Weakly nonlinear case}
\label{sec:wnlgo}

In this section, we prove Proposition~\ref{prop:wnlgo}. 
 Instead of \eqref{eq:NLS}-\eqref{eq:CI}, we consider the weakly
  nonlinear case,
  \begin{equation*}
     i\eps\d_t u^\eps +\frac{\eps^2}{2}\Delta u^\eps=\eps|u^\eps|^2
     u^\eps\quad ;\quad u_0^\eps(x)=\sum_{j=1}^N \alpha_j(x)e^{i\varphi_j(x)/\eps} .
  \end{equation*}
When $d\ge 2$, the creation of new WKB
terms is possible by resonant interactions, provided that $N\ge 3$, as
recalled in the appendix. The
one-dimensional cubic case is special, as there are no nontrivial
resonances, see \cite{CDS10}. In order to present an argument
including the cubic one-dimensional case, we propose a proof which
does not use the creation of, e.g., a fourth term out of three.
\smallbreak

Consider linear phases,  
\begin{equation*}
  \varphi_j(x) = k_j\cdot x.
\end{equation*}
The first part of Proposition~\ref{prop:wnlgo} is simply a restatement
of Theorem~\ref{theo:superp} in this case. To prove the appearance of
nonlinear interferences, we will not consider cutoff functions, and
rely on explicit computations.
WKB analysis in the monokinetic case $N=1$ leads to the hierarchy
\begin{align*}
  & \d_t \phi + \frac{1}{2}|\nabla \phi|^2=0\quad ;\quad \phi
    (0,x)=k\cdot x,\\
  & \d_t a+\nabla \phi\cdot \nabla a+\frac{1}{2}a\Delta \phi =-i|a|^2a\quad
    ;\quad a(0,x)=\alpha(x).
\end{align*}
The  eikonal equation is solved explicitly,
\begin{equation*}
  \phi(t,x) = k\cdot x-\frac{|k|^2}{2}t.
\end{equation*}
As $\Delta \phi=0$, the initial amplitude $\alpha$ is
transported along the vector $k$ with a phase self-modulation:
\begin{equation*}
  a(t,x) = \alpha(x-tk)e^{-it |\alpha(x-tk)|^2}.
\end{equation*}
In the case $N=2$, no new WKB term is created, but interactions between
the two modes lead to a modification of the phase modulation. As
computed in \cite[Section~3]{CDS10}, we find
  \begin{align*}
    a_1(t,x) &= \alpha_1(x-tk_1) e^{-i\(2 \int_0^t
      |\alpha_2(x+(\tau-t)k_1-\tau k_2)|^2  
d\tau + t |\alpha_1(x-tk_1)|^2\)},\\
a_2(t,x) &= \alpha_2(x-tk_2) e^{-i\(2 \int_0^t
  |\alpha_1(x+(\tau-t)k_2-\tau k_1)|^2 d\tau + t
  |\alpha_2(x-t k_2)|^2\)}.  
  \end{align*}
In addition, we have
\begin{equation*}
  \sup_{t\in [0,T]}\left\| u^\eps(t)- a_1(t) e^{i\phi_1(t)/\eps} -
    a_2(t) e^{i\phi_2(t)/\eps} \right\|_{L^2\cap L^\infty}=\O(\eps),
\end{equation*} 
for any $T>0$ (independent of $\eps$); see Corollary~5.13 and
Theorem~6.5 in \cite{CDS10}. 
The leading order nonlinear interactions between the two modes correspond to the
integrals in time in the exponentials. For small time though, the
integrals are zero on the support of the transported amplitudes: in
other words, there exists $T_0>0$ independent of $\eps$ such that for
$t\in [0,T_0]$,
 \begin{align*}
    a_1(t,x) &= \alpha_1(x-tk_1) e^{-i t |\alpha_1(x-tk_1)|^2},\\
a_2(t,x) &= \alpha_2(x-tk_2) e^{-i t
  |\alpha_2(x-t k_2)|^2},
  \end{align*}
in agreement with the conclusion of Theorem~\ref{theo:superp} (up to
the order of precision).

The conclusion in
Proposition~\ref{prop:wnlgo} then follows from the property: 
\begin{equation*}
  \sup_{t\in [0,T_1]} \left\|\alpha_1 (\cdot - tk_1) \(e^{-2i \int_0^t
      |\alpha_2(\cdot+(\tau-t)k_1-\tau k_2)|^2  
d\tau } -
1\)\right\|_{L^p}>0, 
\end{equation*}
or, equivalently,
\begin{equation*}
  \sup_{t\in [0,T_1]} \left\|\alpha_1  \sin \( \int_0^t
      |\alpha_2(\cdot+\tau (k_1- k_2))|^2  
d\tau \)\right\|_{L^p}>0.
\end{equation*}
This is possible as soon as the transport of the support of $\alpha_1$
meets the support of $\alpha_2$, as transported in the above integral. 
Let $\alpha\in C_0^\infty(\R^d)$ supported in the ball centered at
the origin, of radius
$1$, and set
\begin{equation*}
  \alpha_1(x) = \alpha(x),\quad \alpha_2(x) = \alpha(x+3 e_1),
\end{equation*}
where $(e_1,\dots ,e_d)$ is the canonical basis of $\R^d$. Setting
$k_2-k_1= \lambda e_1$ for $\lambda>0$, we see that the above property
is satisfied for some $T_1>0$. We also remark that $T_1\to 0$ as
$\lambda\to \infty$.

\appendix

\section{Weakly nonlinear case and creation of a new term}
\label{sec:wnlgo2}
In this appendix, we prove that in the weakly nonlinear case, if $d\ge
2$, then nonlinear interactions may lead to the creation of a new WKB
term, which is a stronger phenomenon than that used in the proof of
Proposition~\ref{prop:wnlgo}. Consider 
  \begin{equation*}
     i\eps\d_t u^\eps +\frac{\eps^2}{2}\Delta u^\eps=\eps|u^\eps|^2
     u^\eps\quad ;\quad u_0^\eps(x)=\sum_{j=1}^N \alpha_j(x)e^{i\varphi_j(x)/\eps} ,
  \end{equation*}
with now $N=3$, and $d\ge 2$. The
one-dimensional cubic case is special, as there are no nontrivial
resonances, see \cite{CDS10}. Again, we 
consider linear phases,  
\begin{equation*}
  \varphi_j(x) = k_j\cdot x.
\end{equation*}
Recall (\cite{Iturbulent}, see also \cite[Lemma~2.2]{CDS10}) that
the resonant set is defined by
\begin{equation*}
  \text{Res}(n) = \{(j,\ell,m),\quad k_j-k_\ell+k_m=k_n,\
    |k_j|^2-|k_\ell|^2+|k_m|^2=|k_n|^2\}
\end{equation*}
is
characterized as follows: $(k_j,k_\ell,k_m)\in  \text{Res}(n) $  when
the endpoints of the vectors $k_j, k_\ell, k_m, 
k_n$ form four corners of a nondegenerate rectangle
with $k_\ell$ and $k_n$ opposing each other, or when this
quadruplet corresponds to one of the following two degenerate cases:
$(k_j=k_n, k_m=k_\ell)$  
or $(k_j=k_\ell, k_m=k_n)$. Note that we always have
\begin{equation}\label{eq:inclusion}
 \{(j,j,n),\ ((n,j,j),\ a_j\not \equiv 0\} \subset  \text{Res}(n),
\end{equation}
where $a_j$ is the amplitude associated with the phase
\begin{equation*}
  \phi_j(t,x) = k_j\cdot x-\frac{|k_j|^2}{2}t. 
\end{equation*}
\smallbreak

In order for the nonlinearity to create a term associated with a
phase $\phi_4$, out of three phases
associated with wave numbers $k_1$, $k_2$ and $k_3$, we must have 
\begin{equation*}
  k_4: = k_2-k_1+k_3,\quad |k_4|^2 = |k_2|^2-|k_1|^2+|k_3|^2.
\end{equation*}
This resonant condition  is equivalent to the
following conditions:
\begin{equation*}
(k_1-k_2)\cdot (k_1-k_3)=0, 
\end{equation*}
and the endpoints of $k_1$, $k_2$ and $k_3$ are not
aligned (the case of alignment corresponds to the set on the left in
\eqref{eq:inclusion}); this is possible with pairwise different $k_1,k_2,k_3$ and
$k_4\not \in \{k_1,k_2,k_3\}$
provided that $d\ge 2$, see \cite{CDS10} (or
\cite[Section~2.6]{CaBook2}).
For instance if $d=2$, we can choose,
for $\lambda>0$,
\begin{equation*}
  k_1=\lambda(1,1),\quad k_2=\lambda(1,0),\quad
  k_3=\lambda(0,1),\quad\text{hence }k_4=(0,0).  
\end{equation*}
In higher dimension, we simply complete each vector by zero coordinates.
Then a new term, associated with the phase $\phi_4$
may be created by nonlinear resonance. Because of the geometric
characterization of resonances, no other term can be created apart
from this one, since we have completed a rectangle. The creation is
effective only if the associated amplitude does not remain zero. The
equation for the 
corresponding amplitude is 
\begin{equation*}
  \d_t a_4+k_4\cdot \nabla a_4= -i\sum_{(j,\ell,m)\in
    \text{Res(4)}}a_j\bar a_\ell a_m,\quad a_{4\mid t=0}=0.
\end{equation*}
More generally, the term $a_n$ solves
\begin{equation*}
   \d_t a_n+k_n\cdot \nabla a_n= -i\sum_{(j,\ell,m)\in
    \text{Res(n)}}a_j\bar a_\ell a_m,\quad a_{n\mid t=0}=\alpha_n.
\end{equation*}

If we assume that the mode $4$ is not effectively created, that is
$a_4\equiv 0$, then the inclusion \eqref{eq:inclusion} is actually an
equality, and
\begin{equation*}
  \d_t a_j+k_j\cdot \nabla a_j=-2i\sum_{k=1}^3|a_k|^2a_j+i|a_j|^2a_j, \quad j=1,2,3.
\end{equation*}
hence
\begin{equation*}
  a_j(t,x) = \alpha_j(x-tk_j)e^{-iS_j(t,x)},\quad j=1,2,3,
\end{equation*}
for some explicit real-valued phase, whose expression is irrelevant
here (see \cite[Section~3.1]{CDS10} for the formula). 
Given any $T>0$, we may choose $\alpha_1,\alpha_2,\alpha_3$ compactly supported, with
disjoint supports, $k_1,k_2,k_3$ like above, so that
\begin{equation*}
  a_2\bar a_1 a_{3\mid t=T/2}\not \equiv 0. 
\end{equation*}
This shows that the term $a_4$ is actually created, in the sense
that $a_{4}$ does not remain trivial on $[0,T]$.  The error estimate proved in
\cite{CDS12} (see also \cite[Section~2.6]{CaBook2})  yields
\begin{equation*}
  \sup_{t\in [0,T]}\left\|u^\eps(t)-\sum_{j=1}^4
    a_j(t)e^{i\phi_j(t)/\eps}\right\|_{L^2\cap L^\infty} =\O(\eps),
\end{equation*}
hence again the conclusion of Proposition~\ref{prop:wnlgo}. The proof
also implies that for $t\in
[0,T]$, $u^\eps$ has a unique Wigner measure,  given by
\begin{equation*}
  \mu(t,dx,d\xi) = \sum_{j=1}^4
    |a_j(t,x)|^2dx\otimes \delta_{\xi=k_j} .
\end{equation*}

\subsection*{Acknowledgments.} The author wishes to thank Patrick
G\'erard for his remarks which helped improve the consistency of the
paper.

\bibliographystyle{abbrv}
\bibliography{multiphase}

\begin{thebibliography}{10}

\bibitem{ACARMA}
T.~Alazard and R.~Carles.
\newblock Supercritical geometric optics for nonlinear {S}chr\"odinger
  equations.
\newblock {\em Arch. Ration. Mech. Anal.}, 194(1):315--347, 2009.

\bibitem{AlGe07}
S.~Alinhac and P.~G{\'e}rard.
\newblock {\em Pseudo-differential operators and the {N}ash-{M}oser theorem},
  volume~82 of {\em Graduate Studies in Mathematics}.
\newblock American Mathematical Society, Providence, RI, 2007.
\newblock Translated from the 1991 French original by Stephen S. Wilson.

\bibitem{BCLGS-p}
T.~Buckmaster, G.~Cao-Labora, and J.~G{\'o}mez-Serrano.
\newblock Smooth imploding solutions for {3D} compressible fluids.
\newblock Preprint, archived at \url{https://arxiv.org/abs/2208.09445}.

\bibitem{CaBKW}
R.~Carles.
\newblock {WKB} analysis for nonlinear {S}chr\"odinger equations with
  potential.
\newblock {\em Comm. Math. Phys.}, 269(1):195--221, 2007.

\bibitem{CaBook2}
R.~Carles.
\newblock {\em Semi-classical analysis for nonlinear {S}chr\"odinger equations:
  {WKB} analysis, focal points, coherent states}.
\newblock World Scientific Publishing Co. Pte. Ltd., Hackensack, NJ, 2nd
  edition, xiv+352 p. 2021.

\bibitem{CDS10}
R.~Carles, E.~Dumas, and C.~Sparber.
\newblock Multiphase weakly nonlinear geometric optics for {S}chr{\"o}dinger
  equations.
\newblock {\em SIAM J. Math. Anal.}, 42(1):489--518, 2010.

\bibitem{CDS12}
R.~Carles, E.~Dumas, and C.~Sparber.
\newblock Geometric optics and instability for {NLS} and {D}avey-{S}tewartson
  models.
\newblock {\em J. Eur. Math. Soc. (JEMS)}, 14(6):1885--1921, 2012.

\bibitem{CazCourant}
T.~Cazenave.
\newblock {\em Semilinear {S}chr\"odinger equations}, volume~10 of {\em Courant
  Lecture Notes in Mathematics}.
\newblock New York University Courant Institute of Mathematical Sciences, New
  York, 2003.

\bibitem{ChironRousset}
D.~Chiron and F.~Rousset.
\newblock Geometric optics and boundary layers for nonlinear {S}chr\"odinger
  equations.
\newblock {\em Comm. Math. Phys.}, 288(2):503--546, 2009.

\bibitem{Iturbulent}
J.~Colliander, M.~Keel, G.~Staffilani, H.~Takaoka, and T.~Tao.
\newblock Transfer of energy to high frequencies in the cubic defocusing
  nonlinear {S}chr\"odinger equation.
\newblock {\em Invent. Math.}, 181(1):39--113, 2010.

\bibitem{PGX93}
P.~G{\'e}rard.
\newblock Remarques sur l'analyse semi-classique de l'\'equation de
  {S}chr\"odinger non lin\'eaire.
\newblock In {\em S\'eminaire sur les \'Equations aux D\'eriv\'ees Partielles,
  1992--1993}, pages Exp.\ No.\ XIII, 13. \'Ecole Polytech., Palaiseau, 1993.

\bibitem{GMMP}
P.~G{\'e}rard, P.~A. Markowich, N.~J. Mauser, and F.~Poupaud.
\newblock Homogenization limits and {W}igner transforms.
\newblock {\em Comm. Pure Appl. Math.}, 50(4):323--379, 1997.

\bibitem{Gra98}
M.~Grassin.
\newblock Global smooth solutions to {E}uler equations for a perfect gas.
\newblock {\em Indiana Univ. Math. J.}, 47(4):1397--1432, 1998.

\bibitem{Grenier98}
E.~Grenier.
\newblock Semiclassical limit of the nonlinear {S}chr\"{o}dinger equation in
  small time.
\newblock {\em Proc. Amer. Math. Soc.}, 126(2):523--530, 1998.

\bibitem{Casey2018}
C.~Jao.
\newblock Energy-critical {NLS} with potentials of quadratic growth.
\newblock {\em Discrete Contin. Dyn. Syst.}, 38(2):563--587, 2018.

\bibitem{JLM}
S.~Jin, C.~D. Levermore, and D.~W. McLaughlin.
\newblock The semiclassical limit of the defocusing {NLS} hierarchy.
\newblock {\em Comm. Pure Appl. Math.}, 52(5):613--654, 1999.

\bibitem{KVZ09}
R.~Killip, M.~Visan, and X.~Zhang.
\newblock Energy-critical {NLS} with quadratic potentials.
\newblock {\em Comm. Partial Differential Equations}, 34(10-12):1531--1565,
  2009.

\bibitem{LeNgTe18}
N.~Lerner, T.~Nguyen, and B.~Texier.
\newblock The onset of instability in first-order systems.
\newblock {\em J. Eur. Math. Soc. (JEMS)}, 20(6):1303--1373, 2018.

\bibitem{LionsPaul}
P.-L. Lions and T.~Paul.
\newblock Sur les mesures de {W}igner.
\newblock {\em Rev. Mat. Iberoamericana}, 9(3):553--618, 1993.

\bibitem{Majda}
A.~Majda.
\newblock {\em Compressible fluid flow and systems of conservation laws in
  several space variables}, volume~53 of {\em Applied Mathematical Sciences}.
\newblock Springer-Verlag, New York, 1984.

\bibitem{MUK86}
T.~Makino, S.~Ukai, and S.~Kawashima.
\newblock Sur la solution \`a support compact de l'\'{e}quations d'{E}uler
  compressible.
\newblock {\em Japan J. Appl. Math.}, 3(2):249--257, 1986.

\bibitem{MaslovFedoryuk}
V.~P. Maslov and M.~V. Fedoryuk.
\newblock {\em Semi-classical approximation in quantum mechanics. {Transl}.
  from the {Russian} by {J}. {Niederle} and {J}. {Tolar}}, volume~7 of {\em
  Math. Phys. Appl. Math.}
\newblock Kluwer Academic Publishers, Dordrecht, 1981.

\bibitem{MRRS22}
F.~Merle, P.~Rapha\"{e}l, I.~Rodnianski, and J.~Szeftel.
\newblock On blow up for the energy super critical defocusing nonlinear
  {S}chr\"{o}dinger equations.
\newblock {\em Invent. Math.}, 227(1):247--413, 2022.

\bibitem{MR4445442}
F.~Merle, P.~Rapha\"{e}l, I.~Rodnianski, and J.~Szeftel.
\newblock On the implosion of a compressible fluid {I}: {S}mooth self-similar
  inviscid profiles.
\newblock {\em Ann. of Math. (2)}, 196(2):567--778, 2022.

\bibitem{MR4445443}
F.~Merle, P.~Rapha\"{e}l, I.~Rodnianski, and J.~Szeftel.
\newblock On the implosion of a compressible fluid {II}: {S}ingularity
  formation.
\newblock {\em Ann. of Math. (2)}, 196(2):779--889, 2022.

\bibitem{GuyCauchy}
G.~M{\'e}tivier.
\newblock Remarks on the well-posedness of the nonlinear {C}auchy problem.
\newblock In {\em Geometric analysis of PDE and several complex variables},
  volume 368 of {\em Contemp. Math.}, pages 337--356. Amer. Math. Soc.,
  Providence, RI, 2005.

\bibitem{Miller2016}
P.~D. Miller.
\newblock On the generation of dispersive shock waves.
\newblock {\em Phys. D}, 333:66--83, 2016.

\bibitem{Robert}
D.~Robert.
\newblock {\em Autour de l'approximation semi-classique. ({Around}
  semiclassical approximation)}, volume~68 of {\em Prog. Math.}
\newblock Birkh{\"a}user, Cham, 1987.

\bibitem{RV07}
E.~Ryckman and M.~Visan.
\newblock Global well-posedness and scattering for the defocusing
  energy-critical nonlinear {S}chr\"{o}dinger equation in {$\Bbb R^{1+4}$}.
\newblock {\em Amer. J. Math.}, 129(1):1--60, 2007.

\bibitem{ThomannAnalytic}
L.~Thomann.
\newblock Instabilities for supercritical {S}chr\"odinger equations in analytic
  manifolds.
\newblock {\em J. Differential Equations}, 245(1):249--280, 2008.

\end{thebibliography}
\end{document}